\documentclass[twoside, 12pt]{article}
\usepackage{amssymb,mathrsfs,psfrag,graphicx}
\usepackage{CJK,bbm,amsfonts,mathrsfs,amsthm,amssymb,amsmath,indentfirst,color}
\usepackage{hyperref}
\def\v{\varepsilon}

\def\t{\vartheta}

\def\k{\kappa}

\def\g{\gamma}

\def\f{\frac}

\def\di{\displaystyle}
\def\i{\infty}
 \newtheorem{lemma}{\bf Lemma}[section]
       \newtheorem{theorem}[lemma]{\bf Theorem}

       \newtheorem{remark}[lemma]{\bf Remark}

\def\XXint#1#2#3{{\setbox0=\hbox{$#1{#2#3}{\int}$}
\vcenter{\hbox{$#2#3$}}\kern-.5\wd0}}

\def\XXparallel#1#2#3{{\setbox0=\hbox{$#1{#2#3}{\parallel}$}
\vcenter{\hbox{$#2#3$}}\kern-.6\wd0}}

\begin{document}
\date{}
\title{\Large \bf The inviscid limit to a contact discontinuity for the compressible Navier-Stokes-Fourier system  using the relative entropy method}
\author{\small \textbf{Alexis Vasseur},\thanks{The work of A. F. Vasseur was partially supported by the NSF Grant DMS 1209420. E-mail: vasseur@math.utexas.edu.}\quad
   and \textbf{Yi
Wang}\thanks{Y. Wang is supported by NSFC grant No. 11171326 and No. 11322106. E-mail: wangyi@amss.ac.cn.}} \maketitle \small
 
 $^\ast$ Department of Mathematics, University of Texas at Austin, Austin, TX 78712, USA
 
$^{\dag}$  Institute of Applied Mathematics, AMSS,  CAS, Beijing 100190,
China
\\

 {\bf Abstract:} We consider the zero heat conductivity  limit to a contact discontinuity for the mono-dimensional  full compressible
Navier-Stokes-Fourier  system. The method is based on the  relative entropy method, and do not assume any smallness conditions on the discontinuity, nor on the $BV$ norm of the initial data. It is proved that for any viscosity $\nu\geq0$, the solution of the compressible
Navier-Stokes-Fourier system (with well prepared initial value) converges,  when the heat
conductivity $\k$ tends to zero,  to the contact discontinuity solution to the corresponding Euler system.  We obtain the  decay rate $\k^{\frac12}$. It implies that the heat conductivity dominates the dissipation in the regime of the limit to a  contact discontinuity. This is the first  result, based on the relative entropy, of an   asymptotic limit to a discontinuous solutions for a system.

\renewcommand{\theequation}{\arabic{section}.\arabic{equation}}
\setcounter{equation}{0}
\section{Introduction}
 We study the zero dissipation limit of the solution to the Navier-Stokes-Fourier
 system of a compressible heat-conducting gas in the Lagrangian coordinates in one-dimensional space $\mathbb{R}$:
\begin{equation}
 \left\{
  \begin{array}{l}
    \di \tau_t-v_x=0,\\
    \di v_t+p(\tau,\vartheta)_x=\nu(\f{v_x}{\tau})_x,\\
    \di (e(\tau,\t)+\f{v^2}{2})_t+(p(\tau,\vartheta)v)_x=\k(\f{\t_x}{\tau})_x+\nu(\f{vv_x}{\tau})_x
  \end{array}
\right. \label{NS}
\end{equation}
where the functions
$\tau(x,t)>0,v(x,t),\t(x,t)>0$ represent the specific volume, velocity
and the absolute temperature of the gas, respectively, 
$p=p(\tau,\t)$, $e=e(\tau,\t)$ is the pressure and the internal energy satisfying the second law of thermodynamics,
$\nu>0$ is the viscosity constant and $\k>0$ is the coefficient of
heat conduction.

If $\k=\nu=0$ in \eqref{NS}, then the corresponding Euler system reads as
\begin{equation}
\left\{
\begin{array}{l}
\tau_t-v_x=0,\\
v_t+p(\tau,\t)_x=0,\\[2mm]
\displaystyle{\Big( e(\tau,\t)+\f{v^2}{2}\Big)_t+(pv)_x=0.}
\end{array}
\right.\label{euler}
\end{equation}

For the Euler equations \eqref{euler}, it is well-known that there are three basic
wave patterns, i.e.,  shock and rarefaction waves in the genuinely nonlinear fields and contact discontinuities in the linearly degenerate field.
If we consider the Riemann problem of compressible Euler system \eqref{euler} with the following Riemann initial data
\begin{equation}\label{RI}
(\bar\tau,\bar v,\bar \t)(x,0)=\left\{
\begin{array}{ll}
\di (\tau_L,v_L,\t_L)&\di x<0,\\
\di (\tau_R,v_R,\t_R)& \di x>0,
\end{array}
\right.
\end{equation}
then the above Riemann problem admits a contact discontinuity solution
\begin{equation}\label{CD}
(\bar\tau,\bar v,\bar \t)(x,t)=\left\{
\begin{array}{ll}
\di (\tau_L, v_L,\t_L)&\di x<\bar v t,\\
\di (\tau_R,v_R,\t_R)& \di x>\bar vt,
\end{array}
\right.
\end{equation}
provided that 
$$
v_L=v_R:=\bar v,\qquad p_L=p_R:=\bar p.
$$

The relative entropy method was first introduced by Dafermos \cite{Da} and Diperna \cite{Di} to prove the $L^2$ stability and uniqueness of Lipschitzian solutions to the conservation laws. The method was used by  Leger in \cite{L} to show the stability of shocks in $L^2$ (up to a shift) in the scalar case (see also \cite{A} for an extension  to $L^p$, $1<p<\infty$).  The method has been extended to the system case in  \cite{LV} for extreme shocks. The method is purely nonlinear, and allow to work without  smallness assumptions.  The relative entropy method  is also an effective method for the study of  asymptotic limits. One of the first use of the method in this context is due to  Yau \cite{Yau} for the hydrodynamic limit of Ginzburg-Landau models. Since then, there have been many works  in this context, see \cite{BGL1, BGL2, BTV, BV, GS, LM, MS} etc. and the survey paper \cite{V}, although they are all concerning the limit to a smooth (Lipschitz) limit function. Recently, the relative entropy method has been successfully applied in  to prove the vanishing viscosity limit of the viscous scalar conservation laws to shocks \cite{CV}, and the $L^2-$contraction of viscous shock profiles to the scalar viscous conservation laws \cite{KaV}. The theory of stability of discontinuous solutions, based on the relative entropy has been reformulated in \cite{SV, KaV2} in terms of contraction, up to a shift. 

\vskip0.1cm
This paper is the first attempt to use this theory for the study of asymptotic limits to a discontinuous solution, in the case  of systems.  In the present paper, we are concerned with the vanishing physical viscosity limit of compressible Navier-Stokes-Fourier system \eqref{NS} to the above contact discontinuity solution \eqref{CD} to the Euler system \eqref{euler}, with arbitrarily large wave strength, using the relative entropy method.  The contact discontinuities of the Euler system enjoys a contraction property of a special inhomogeneous pseudo-distance, which does not need a shift (see \cite{SV2}). This property is one of the key of our proof.  

For the zero dissipation limit of conservation laws, there are many results using BV-theory, see Bianchini-Bressan\cite{BB} etc. and the limit to basic wave patterns to the compressible Navier-Stokes equations and Boltzmann equation by the elementary energy method, one can refer to \cite{GX, M, W, HWY, HWY1, HWWY}  and the references therein. 
\vskip0.1cm

In the present paper, we consider the state equation as in Feireisl \cite{F}
\begin{equation}\label{pe}
p(\tau,\t)=p_e(\tau)+p_{th}(\tau,\t)
\end{equation}
with $p_e(\tau)$ being the elastic pressure. The thermal pressure $p_{th}(\tau,\t)$   takes the forms as
$$
p_{th}(\tau,\t)=\t p_\t(\tau).
$$
Both functions $p_e(\tau)$ and $p_\t(\tau)$  satisfy
\begin{equation}\label{c1}
p_e(\tau), \ p_\t(\tau)\geq0,\qquad p_e^\prime(\tau),\  p_\t^\prime(\tau)\leq 0.
\end{equation}
 Correspondingly, the internal energy takes the form as
 \begin{equation}\label{ee}
 e(\tau,\t)=P_e(\tau)+Q(\t),
 \end{equation}
 where $P_e(\tau)$ is the elastic potential defined by
 $$
 P_e(\tau)=-\int^\tau p_e(y)dy,
 $$
 and $Q(\t)=\int^\t C_v(z)dz$ is the thermal energy satisfying
 $$
 \frac{\partial Q}{\partial \t}=e_\t=C_v(\t),
 $$
 where $C_v(\t)$ is a positive and smooth function of $\t.$ 
 By the second law of thermodynamic
  $$
  \t ds=de+pd\tau,
  $$
  we can calculate the entropy $s$ as
\begin{equation}\label{entropy}
s(\tau,\t)=\int^\t\frac{C_v(z)}{z}dz-P_\t(\tau),
\end{equation}
with 
$$
P_\t(\tau)=-\int^\tau p_\t(\tau)d\tau.
$$

Set $U=(\tau,v,e(\tau,\t)+\f{v^2}{2})^t$ being the solution to the compressible Navier-Stokes equations \eqref{NS} and $\bar U=(\bar\tau,\bar v,\bar e(\bar\tau,\bar\t)+\f{\bar v^2}{2})^t$ being the contact discontinuity \eqref{CD}  to the compressible Euler equations \eqref{euler}. As usual, the associated relative entropy is defined by
$$
s(U|\bar U):=s(U)-s(\bar U)-ds(\bar U)\cdot (U-\bar U),
$$
with $s$ being the entropy in \eqref{entropy}.  Correspondingly, the relative entropy-flux is 
$$
q(U;\bar U):=-ds(\bar U)\cdot (f(U)-f(\bar U)),
$$
with the flux $f(U)=(-v,p(\tau,\t), p(\tau,\t)v)^t.$

Define the entropy functional as
\begin{equation}\label{ef}
\mathcal{E}(t)=-\int\bar\t\ s(U|\bar U) dx.
\end{equation}
Direct calculations show that for the state equations in \eqref{pe} and \eqref{ee}, one has the explicit formula for \eqref{ef}
\begin{equation}\label{s1}
\begin{array}{ll}
\di \mathcal{E}(t)=-\int\bar\t~s(U|\bar U)dx=\int\big[-\bar\t~s(U)+\bar \t~s(\bar U)+\bar \t~ds|_{U=\bar U}\cdot(U-\bar U)\big]dx\\
\di =\int \Big[\bar \t P_\t(\tau)-\bar\t\int^\t \f{C_v(z)}{z}dz-\bar\t P_\t(\bar\tau)+\bar\t\int^{\bar\t} \f{C_v(z)}{z}dz\\
\di\qquad\qquad +\bar p(\tau-\bar\tau)-\bar v(v-\bar v)+(E-\bar E)\Big]dx\\
\di =\int\Big[\bar\t\big(P_\t(\tau)-P_\t(\bar\tau)+p_\t(\bar\tau)(\tau-\bar\tau)\big)+\big(P_e(\tau)-P_e(\bar\tau)+p_e(\bar\tau)(\tau-\bar\tau)\big)\\
\di\qquad\qquad +\big(Q(\t)-Q(\bar\t)-\bar\t\int^{\t}_{\bar\t}\f{C_v(z)}{z}dz\big)+\f12(v-\bar v)^2\Big]dx.
\end{array}
\end{equation}

Now we state the main result as follows.
%%%%%%%%%%%%%%%%%%%%%%%%%%%%%%%%%%%%%%%%%%%%%%%%%%%%%%%%%%%%%%%%%%%%%%%%%%%%%%
\begin{theorem}\label{limit-th} For $\k>0$, let $(\tau_\k, v_\k,\t_\k)$ be a weak solution to \eqref{NS} such that $s(U_\k|\bar U)$ lies in $L^1(0,T; L^1(\mathbf{R}))$  where $U_\k= (\tau_\k, v_\k,\t_\k)$ and $\bar U$ is the contact discontinuity to the Euler system \eqref{euler} and $s$ is defined in \eqref{entropy}.  Assume further that the pressure $p$ satisfies that
\begin{equation}\label{pe-c}
p_e(\tau)\sim \tau^{-\gamma},\quad {\rm with}~ \gamma\geq2, \quad{\rm as}~\tau\rightarrow 0+,
\end{equation}
or
\begin{equation}\label{pt-c}
p_\t(\tau)\sim \tau^{-\gamma},\quad {\rm with}~ \gamma\geq2, \quad{\rm as}~\tau\rightarrow 0+,
\end{equation}
Then if the initial perturbation satisfies that 
$$
\mathcal{E}(0)=O(\sqrt{\kappa}),
$$
and
\begin{equation}\label{ic}
s(U_{0\k}|\bar U_0)\in L^\i(\Omega)
\end{equation}
with $\Omega\subset \mathbf{R}$ being any neighborhood of $x=0$. Then for any $\nu\geq0$, there exists $\kappa_0>0$ such that if $0<\kappa\leq \kappa_0$, then it holds that
$$
\mathcal{E}(t)=O(\sqrt{\kappa}),\quad  \forall t\in[0,T]~{\it with~ any~ fixed}~ T>0.
$$
\end{theorem}

\begin{remark} 
There is no need of uniform-in-$\k$ $L^1$ bound on the weak solution $U_\k$, even though its existence is still open now. However, if both the initial volume and initial temperature have the lower and upper bound, then existence of global classical solution can be proved by following the similar procedure as in Kazhikhov-Shelukhin \cite{KS}.
\end{remark}

\begin{remark} 
By \eqref{s1}, it can be seen that $\|v_\k-\bar v\|^2=O(\sqrt\k)$. However, for $\tau_\k-\bar\tau$ and $\t_\k-\bar\t$, one only can get the convergence in entropy norm as in \eqref{s1} and  can not derive the convergence in $L^2-$norm unless one has the uniform lower and upper bound for both $\tau_\k$ and $\t_\k$ a priorly.
\end{remark}

\begin{remark} Theorem \ref{limit-th} shows that the solution of compressible Navier-Stokes system can converge to contact discontinuity solution to the corresponding Euler system with the decay rate $\sqrt{\k}$, uniformly for any viscosity $\nu\geq0$, which implies something surprising that the heat conductivity $\kappa$ dominate the dissipation in the regime of limit of contact discontinuity as $\k$ tends to zero. Note that the wave strength of contact discontinuity can be arbitrarily large and  the initial perturbation can contain the possible vacuum states formally.
\end{remark}

\begin{remark}
For the state equation, the pressure $p$ can contain a large class of gases, such as
$$
p(\tau,\t)=\f{R\t}{\tau}+p_e(\tau)
$$
with $R>0$ being the gas constant and the elastic pressure
$p_e(\tau)$ satisfying the condition \eqref{pe-c}.
\end{remark}

\section{Proof of Theorem 1.1}
\setcounter{equation}{0}
In this section, we will prove Theorem 1.1 by relative entropy method. We will first define the cut-off entropy functional by suitably choosing the cut-off function $\eta(x)$ based on the underlying structure of contact discontinuity.
The cut-off function $\eta(x)$ is chosen to satisfy
\begin{equation}\label{eta}
\eta(x)=\left\{
\begin{array}{ll}
\di 1,&\di x\leq -1,\\
\di {\rm smooth~ and ~decreasing},&\di -1\leq x\leq1,\\
\di 0,&\di x\geq1,
\end{array}
\right.
\end{equation}
and
\begin{equation}\label{eta-p}
\eta^\prime(-x)=\eta^\prime (x).
\end{equation}
For definiteness, we can choose $\eta(x)$ as follows.
\begin{equation}
\eta(x)=\left\{
\begin{array}{ll}
\di 1,&\di x\leq -1,\\
\di \f14 x^3-\f34x+\f12,&\di -1\leq x\leq1,\\
\di 0,&\di x\geq1.
\end{array}
\right.
\end{equation}
Obviously, we have
\begin{equation}\label{eta-1}
\eta^\prime(x)=\left\{
\begin{array}{ll}
\di \f34x^2-\f34,&\di -1\leq x\leq1,\\
\di 0,&\di {\rm otherwise},
\end{array}
\right.
\end{equation}
satisfying $\eta^\prime\leq0,~{\rm supp} ~\eta^\prime=[-1,1]$ and  the symmetry in \eqref{eta-p}.

Define the cut-off entropy functional as
\begin{equation}\label{E1}
\mathcal{E}_1(t)=-\int_{-\infty}^{+\infty}\eta(\f{x}{\v})\t_L~ s(U|U_L)dx-\int_{-\infty}^{+\infty}\eta(-\f{x}{\v})\t_R ~s(U|U_R)dx
\end{equation}
where $\v>0$ represents the layer-width for the contact discontinuity in the vanishing dissipation limit and is to be determined as $\v=\sqrt\kappa$ in the sequel. The functional $\mathcal{E}_1(t)$ and the cut-off function $\eta(x)$ both play an important role in our analysis.

In fact, the following proof is independent of the state equations in \eqref{pe} and \eqref{ee} until the estimation \eqref{E131}, that is, all the proofs before \eqref{E131} are valid to the generic gas satisfying the second law of thermodynamic. By the second law of thermodynamics, it holds that
$$
\t ds=de+pd\tau=d(E-\f{v^2}{2})+pd\tau=dE-vdv+pd\tau,
$$
where the total energy $E=e(\tau,\t)+\f{v^2}{2}$, thus the entropy $s$ satisfies the equation
\begin{equation}\label{s-e}
s_t=\k(\f{\t_x}{\tau\t})_x+\k\f{\t_x^2}{\tau\t^2}+\nu \f{v_x^2}{\tau\t}.
\end{equation}

By direct calculations, it holds that
\begin{equation}\label{E1-1}
\begin{array}{ll}
\di \f{d\mathcal{E}_1(t)}{dt}= -\int_{-\infty}^{+\infty}\bigg(\eta(\f x\v)\t_L+\eta(-\f x\v)\t_R\bigg)\bigg(\nu\f{v_x^2}{\tau\t}+\k\f{\t_x^2}{\tau\t^2}\bigg)dx\\
\di -\int_{-\infty}^{+\infty}\bigg(\eta(\f x\v)\t_L+\eta(-\f x\v)\t_R\bigg)\k\big(\f{\t_x}{\tau\t}\big)_xdx\\
\di -\int_{-\infty}^{+\infty}\Big(\eta(\f x\v)+\eta(-\f x\v)\Big)\Big(\bar pv_x-\bar v\big(\nu(\f{v_x}{\tau})_x-p_x\big)+\k\big(\f{\t_x}{\tau}\big)_x+\nu\big(\f{vv_x}{\tau}\big)_x-(pv)_x\Big)dx.
\end{array}
\end{equation}
Note that the third term on the right hand side of \eqref{E1-1} vanishes due to the integration by part and the symmetry of the $\eta^\prime(x)$. In fact, 
\begin{equation}\label{E3}
\begin{array}{ll}
\di-\int_{-\infty}^{+\infty}\Big(\eta(\f x\v)+\eta(-\f x\v)\Big)\Big(\bar pv_x-\bar v\big(\nu(\f{v_x}{\tau})_x-p_x\big)+\k\big(\f{\t_x}{\tau}\big)_x+\nu\big(\f{vv_x}{\tau}\big)_x-(pv)_x\Big)dx\\
\di =-\int_{-\infty}^{+\infty}\Big(\eta(\f x\v)+\eta(-\f x\v)\Big)\Big(\bar p(v-\bar v)_x-\bar v\big(\nu(\f{v_x}{\tau})_x-(p-\bar p)_x\big)\\
\di \qquad\qquad\qquad\qquad\qquad\qquad\qquad +\k\big(\f{\t_x}{\tau}\big)_x+\nu\big(\f{vv_x}{\tau}\big)_x-(pv-\bar p\bar v)_x\Big)dx\\
\di\di =-\Big(\eta(\f x\v)+\eta(-\f x\v)\Big)\Big(\bar p(v-\bar v)-\bar v\nu\f{v_x}{\tau}-(p-\bar p)+\k\f{\t_x}{\tau}+\nu\f{vv_x}{\tau}-(pv-\bar p\bar v)\Big)|_{-\infty}^{+\infty}\\
\di+\int_{-\infty}^{+\infty}\f1\v\Big(\eta^\prime(\f x\v)-\eta^\prime(-\f x\v)\Big)\Big(\bar p(v-\bar v)-\bar v\nu\f{v_x}{\tau}-(p-\bar p)+\k\f{\t_x}{\tau}+\nu\f{vv_x}{\tau}-(pv-\bar p\bar v)\Big)dx\\
\di \equiv0.
\end{array}
\end{equation}
By integrating \eqref{E1-1} with respect to $t$ over $[0,t]$ with $t\in[0,T]$,  wow we arrive at 
\begin{equation}\label{E1-2}
\begin{array}{ll}
\di \mathcal{E}_1(t) +\int_{-\infty}^{+\infty}\bigg(\eta(\f x\v)\t_L+\eta(-\f x\v)\t_R\bigg)\bigg(\nu\f{v_x^2}{\tau\t}+\k\f{\t_x^2}{\tau\t^2}\bigg)dx\\
\di =\mathcal{E}_1(0)-\int_0^t\int_{-\infty}^{+\infty}\bigg(\eta(\f x\v)\t_L+\eta(-\f x\v)\t_R\bigg)\k\big(\f{\t_x}{\tau\t}\big)_xdxdt.\\
\end{array}
\end{equation}
Note that the terms on the right hand side of \eqref{s1}  are non-negative under the condition \eqref{c1}.
Then it holds that 
\begin{equation}\label{s3}
\begin{array}{ll}
\di \mathcal{E}_1(t)=-\int_{-\infty}^{+\infty}\eta(\f{x}{\v})\t_L~ s(U|U_L)dx-\int_{-\infty}^{+\infty}\eta(-\f{x}{\v})\t_R ~s(U|U_R)dx\\
\di \geq \int_{-\infty}^{0}\eta(\f{x}{\v})[-\t_L~ s(U|U_L)]dx+\int_{0}^{+\infty}\eta(-\f{x}{\v})[-\t_R ~s(U|U_R)]dx\\
\di\geq \f12 \int_{-\infty}^{0}[-\t_L~ s(U|U_L)]dx+\f12\int_{0}^{+\infty}[-\t_R ~s(U|U_R)]dx\\
\di =\f12\int_{-\infty}^{+\infty}[-\bar\t ~s(U|\bar U)]dx=\f12\mathcal{E}(t).
\end{array}
\end{equation}
Substituting the above inequalities \eqref{s1}-\eqref{s3} into \eqref{E1-2}, one has
\begin{equation}\label{E1-21}
\begin{array}{ll}
\di \f12 \int_{-\infty}^{+\infty}\Big[ \bar \t\big(P_\t(\tau)-P_\t(\bar \tau)+p_\t(\bar\tau)(\tau-\bar\tau)\big)+\big(P_e(\tau)-P_e(\bar \tau)+p_e(\bar\tau)(\tau-\bar\tau)\big)\\
\di\qquad +\big(Q(\t)-Q(\bar\t)-\bar\t\int^{\t}_{\bar\t}\f{C_v(z)}{z}dz\big)+\f12(v-\bar v)^2\Big](x,t)dx\\
\di\qquad +\f12\min\{\t_L,\t_R\}\int_0^t\int_{-\infty}^{+\infty}\big(\nu\f{v_x^2}{\tau\t}+\k\f{\t_x^2}{\tau\t^2}\big)dxdt\\
\di \leq \mathcal{E}_1(0)-\int_0^t\int_{-\infty}^{+\infty}\bigg(\eta(\f x\v)\t_L+\eta(-\f x\v)\t_R\bigg)\k\big(\f{\t_x}{\tau\t}\big)_xdxdt.
\end{array}
\end{equation}
Now we treat the terms on the right hand side of \eqref{E1-21}. First, for the initial values satisfying \eqref{ic}, it holds that
\begin{equation}\label{E0}
\begin{array}{ll}
\di \mathcal{E}_1(0)=-\int_{-\infty}^{+\infty}\eta(\f{x}{\v})\t_L~ s(U_0|U_L)dx-\int_{-\infty}^{+\infty}\eta(-\f{x}{\v})\t_R ~s(U_0|U_R)dx\\
\di =\int_{-\infty}^{0}\eta(\f{x}{\v})[-\t_L~ s(U_0|U_L)]dx+\int_{0}^{\v}\eta(\f{x}{\v})[-\t_L~ s(U_0|U_L)]dx\\
\di \quad +\int_{0}^{+\infty}\eta(-\f{x}{\v})[-\t_R ~s(U_0|U_R)]dx+\int_{-\v}^{0}\eta(-\f{x}{\v})[-\t_R ~s(U_0|U_R)]dx\\
\di \leq \int_{-\infty}^{+\infty}[-\bar \t_0~ s(U_0|\bar U_0)]dx+\int_{0}^{\v}\eta(\f{x}{\v})[-\t_L~ s(U_0|U_L)]dx\\
\di\qquad\qquad +\int_{-\v}^{0}\eta(-\f{x}{\v})[-\t_R ~s(U_0|U_R)]dx\\
\di \leq \mathcal{E}(0)+C\int_0^\v\eta(\f x\v)dx+C\int_{-\v}^0\eta(-\f x\v)dx \leq \mathcal{E}(0)+C\v
\end{array}
\end{equation}
for some generic positive constant $C$ provided that $\v=\sqrt\k$ is suitably small.
Now it remains to control the last term on the right hand side of \eqref{E1-21}. In fact, by integration by part and Cauchy inequality, it holds that
\begin{equation}\label{E1-3}
\begin{array}{ll}
\di -\int_0^t\int_{-\infty}^{+\infty}\bigg(\eta(\f x\v)\t_L+\eta(-\f x\v)\t_R\bigg)\k\big(\f{\t_x}{\tau\t}\big)_xdxdt\\
\di =\int_0^t\int_{-\infty}^{+\infty}\bigg(\eta^\prime(\f x\v)\t_L-\eta^\prime(-\f x\v)\t_R\bigg)\f{\k}{\v}\f{\t_x}{\tau\t}dxdt=\int_0^t\int_{-\v}^{\v}\eta^\prime(\f x\v)\big(\t_L-\t_R\big)\f{\k}{\v}\f{\t_x}{\tau\t}dxdt\\
\di \leq\f14 \min\{\t_L,\t_R\} \int_0^t\int_{-\v}^{\v}\k\f{\t_x^2}{\tau\t^2}dxdt+\f{(\t_L-\t_R)^2}{\min\{\t_L,\t_R\}}\frac{\k}{\v^2}\int_0^t\int_{-\v}^{\v}\big(\eta^\prime(\f x\v)\big)^2\f1\tau dxdt\\
\di \leq\f14 \min\{\t_L,\t_R\} \int_0^t\int_{-\v}^{\v}\k\f{\t_x^2}{\tau\t^2}dxdt+C(\t_L,\t_R)\frac{\k}{\v^2}\int_0^t\int_{-\v}^{\v}\f1\tau dxdt
\end{array}
\end{equation}
where the positive constant $C(\t_L,\t_R)=\f{9(\t_L-\t_R)^2}{16\min\{\t_L,\t_R\}}.$
With the notation ${\bf 1}_{\{\cdots\}}$ being the characteristic function of the set $\{\cdots\}\subset\mathbf{R}\times[0,T]$ and $C_1$ being suitably large positive constant  to be determined, it holds that
\begin{equation}\label{E131}
\begin{array}{ll}
\di \int_0^t\int_{-\v}^{\v}\f1\tau dxdt=\int_0^t\int_{-\v}^{\v}\f1\tau\big({\bf 1}_{\{\f1\tau\leq C_1\}}+{\bf 1}_{\{\f1\tau> C_1\}}\big) dxdt\\
\di\leq 2C_1\v+\int_0^t\int_{-\v}^{\v}\f1\tau {\bf 1}_{\{\f1\tau> C_1\}}dxdt.
\end{array}
\end{equation}
Now the state equations in \eqref{pe} and \eqref{ee} as in \cite{F} are crucially used to control the second term on the right hand of \eqref{E131}.  Moreover, if the pressure function $p_e(\tau)$ satisfies \eqref{pe-c} or $p_\t(\tau)$ satisfies \eqref{pt-c}, then
$$P_e(\tau)\ 
{\rm or}\  P_\t(\tau)\sim \tau^{-(\g-1)},\quad \g\geq2,\quad {\rm as}~ \tau\rightarrow0+.
$$
Therefore, we can choose the suitably large positive constant $C_1$ such that
\begin{equation}\label{p-c}
|\frac {\frac1\tau}{\bar \t\big(P_\t(\tau)-P_\t(\bar \tau)+p_\t(\bar\tau)(\tau-\bar\tau)\big)+\big(P_e(\tau)-P_e(\bar \tau)+p_e(\bar\tau)(\tau-\bar\tau)\big)}| {\bf 1}_{\{\f1\tau> C_1\}}<+\infty,
\end{equation}
therefore, one has
\begin{equation}\label{re}
\begin{array}{ll}
\di |\int_0^t\int_{-\v}^{\v}\f1\tau {\bf 1}_{\{\f1\tau> C_1\}}dxdt|\\
\di\leq C\int_0^t\int_{-\i}^{\i}\big[\bar \t\big(P_\t(\tau)-P_\t(\bar \tau)+p_\t(\bar\tau)(\tau-\bar\tau)\big)\\
\di\qquad\qquad\qquad\qquad\qquad  +\big(P_e(\tau)-P_e(\bar \tau)+p_e(\bar\tau)(\tau-\bar\tau)\big)\big]dxdt.
\end{array}
\end{equation}
Substituting \eqref{re} into \eqref{E131} and then into \eqref{E1-3} yields that
\begin{equation}\label{E1-f}
\begin{array}{ll}
\di -\int_0^t\int_{-\infty}^{+\infty}\bigg(\eta(\f x\v)\t_L+\eta(-\f x\v)\t_R\bigg)\k\big(\f{\t_x}{\tau\t}\big)_xdxdt\\
\di \leq\f14 \min\{\t_L,\t_R\} \int_0^t\int_{-\v}^{\v}\k\f{\t_x^2}{\tau\t^2}dxdt +C(\t_L,\t_R)\frac{\k}{\v^2}\int_0^t\int_{-\i}^{\i}\big[\bar \t\big(P_\t(\tau)-P_\t(\bar \tau)\\[3mm]
\di\qquad\qquad\qquad\qquad\qquad\qquad  +p_\t(\bar\tau)(\tau-\bar\tau)\big)+\big(P_e(\tau)-P_e(\bar \tau)+p_e(\bar\tau)(\tau-\bar\tau)\big)\big]dxdt.
\end{array}
\end{equation}
The combination of  \eqref{E1-21}, \eqref{E0}  and \eqref{E1-f} implies that
\begin{equation}\label{f}
\begin{array}{ll}
\di \f12\mathcal{E}(t)+\f14\min\{\t_L,\t_R\}\int_0^t\int_{-\infty}^{+\infty}\big(\nu\f{v_x^2}{\tau\t}+\k\f{\t_x^2}{\tau\t^2}\big)dxdt\\
\di\qquad\qquad\qquad\qquad \leq \mathcal{E}(0)+C\v+C\frac{\k}{\v^2}\Big(\v+\int_0^t\mathcal{E}(t)dt\Big).
\end{array}
\end{equation}
Set
$$
F(t)=\int_0^t\mathcal{E}(t)dt,
$$
then applying Gronwall inequality to \eqref{f}  implies that 
\begin{equation}
F(t)\leq C_2\Big(\v+\f{\mathcal{E}(0)+\v}{\f\k{\v^2}}\Big)e^{C\f\k{\v^2}t},
\end{equation}
for some positive constant $C_2$.
Furthermore, one has
\begin{equation}
\mathcal{E}(t)\leq C\Big[\mathcal{E}(0)+\v+\frac{\k}{\v^2}\Big(\v+\Big(\v+\f{\mathcal{E}(0)+\v}{\f\k{\v^2}}\Big)e^{C\f\k{\v^2}t}\Big)\Big].
\end{equation}
Choosing the contact discontinuity layer $\v=\sqrt{\k}$ and the initial perturbation $\mathcal{E}(0)=O(\sqrt{\k})$ as in Theorem 1.1, one has
$$
\mathcal{E}(t)=O(\sqrt{\k}), \qquad \forall t\in[0,T], ~~\forall \nu\geq0,
$$
which completed the proof of Theorem 1.1.
%\vspace{4mm}
%\noindent {\bf Acknowledgment:}\,\,

\end{document}